\magnification =\magstep1
\baselineskip=13.0pt
\overfullrule=0pt

\font\sc=cmcsc10 at 12pt
\font\large=cmbx12 at 14pt

\def\Ass{\hbox{Ass}\hskip -0.0em}
\def\qed{\hbox{\quad \vrule width 1.6mm height 1.6mm depth 0mm}\vskip 1ex}
\def\eqed{\hbox{\quad \vrule width 1.6mm height 1.6mm depth 0mm}}
\def\mvrule{\vrule height 3.2ex depth 1.6ex width .01em}
\def\cl{:\quad}

\def \thmline#1{\vskip 20pt
	\noindent{\sc #1:}\ \ %
	\bgroup \advance\baselineskip by -1pt \it
	\abovedisplayskip =4pt
	\belowdisplayskip =3pt
	\parskip=0pt
	}
\def \endb{\egroup \vskip 1.4ex}
\def \proof{\smallskip\noindent {\sl Proof:\ \ }}

\ %
\vskip 3ex
\leftline{\large THE FIRST MAYR-MEYER IDEAL}
\vskip 2ex
\noindent{\sc IRENA SWANSON}
\vskip 2ex

\bgroup
\baselineskip=11pt
\noindent New Mexico State University - Department of Mathematical Sciences,
Las Cruces, New Mexico 88003-8001, USA.
E-mail: {\tt iswanson@nmsu.edu}.

\egroup

\vskip 5ex
{\bgroup \narrower \narrower
\baselineskip=11pt
\noindent
{\bf Summary.}
This paper gives a complete primary decomposition
of the first, that is, the smallest, Mayr-Meyer ideal,
its radical, and the intersection of its minimal components.
The particular membership problem which makes the Mayr-Meyer ideals'
complexity doubly exponential in the number of variables
is here examined also for the radical and the intersection of the
minimal components.
It is proved that for the first Mayr-Meyer ideal the complexity
of this membership problem is the same as for its radical.
This problem was motivated by a question of Bayer, Huneke and Stillman.

}

\vskip 4ex

Grete Hermann proved in [H] that for any ideal $I$
in an $n$-dimensional polynomial ring over the field of rational numbers,
if $I$ is generated by polynomials $f_1, \ldots, f_k$ of degree at most $d$,
then it is possible to write $f = \sum r_i f_i$,
where each $r_i$ has degree at most $\deg f + (kd)^{(2^n)}$.
Mayr and Meyer in [MM] found ideals $J(n,d)$
for which a doubly exponential bound in $n$ is indeed achieved.
Bayer and Stillman [BS] showed that for these same ideals
also any minimal generating set of syzygies
has elements of degree which is doubly exponential in $n$.
Koh [K] modified the original ideal to obtain homogeneous quadric ideals
with doubly exponential degrees of syzygies and ideal membership coefficients.

Bayer, Huneke and Stillman have raised questions about the structure
of these Mayr-Meyer ideals:
is the doubly exponential behavior due to the number of minimal primes,
to the number of associated primes,
or to the structure of one of them?
This paper, together with [S],
is an attempt at answering these questions.
More precisely,
the Mayr-Meyer ideal $J(n,d)$ is an ideal in a polynomial ring in $10n + 2$
variables whose generators have degree at most $d+2$.
This paper analyzes the case $n = 1$
and shows that in this base case the embedded components do not play a role.

Theorem~1 of this paper gives a complete primary decomposition
of $J(1,d)$,
after which the intersection of the minimal components
and the radical come as easy corollaries.
The last proposition shows that the complexity of the particular membership
problem from [MM, BS, K] for the radical of $J(1,d)$
is the same as the complexity of the membership problem for $J(1,d)$.
Thus at least for the case $n = 1$,
neither the embedded components nor the non-reducedness play a role
in the complexity.

In a developing paper ``Primary decomposition of the Mayr-Meyer ideal" [S],
partial primary decompositions
are determined for the Mayr-Meyer ideals $J(n,d)$ for all $n \ge 2$, $d \ge 1$.
Under the assumption that the characteristic of the field does not divide $d$,
for $n \ge 2$,
the number of minimal primes is exactly $nd^2 + 20$,
and the number of embedded primes likewise depends on $n$ and $d$.
However,
a precise number of embedded components is not known.
The case $n = 1$ is very different from the case $n \ge 2$.
For example,
under the same assumption on the characteristic of the field,
the number of minimal primes of the first Mayr-Meyer ideal is $d+4$,
and there is exactly one embedded prime.
For understanding the asymptotic behavior of the Mayr-Meyer ideals $J(n,d)$,
the case $n = 1$ may not seem interesting,
however, it is a basis of the induction arguments
for the behavior of the other $J(n,d)$.
Furthermore,
the case $n = 1$ is computationally and notationally more accessible.

All results of this paper were verified for specific low values of $d$
on Macaulay2.

{\bf Acknowledgement.}
This research was done during 2000/01 at University of Kansas,
supported by the NSF-POWRE grant.
I thank Craig Huneke for introducing me to the Mayr-Meyer ideals,
for all the conversations about them,
and for his enthusiasm for this project.
I also thank the NSF for partial support on DMS-9970566.

\vskip 3ex
The first Mayr-Meyer ideal $J(1,d)$ is defined as follows.
Let $K$ be a field, and $d$ a positive integer.
In case the characteristic of $K$ is a positive prime $p$,
write $d = d'i$, where $i$ is a power of $p$,
and $d'$ and $p$ are relatively prime integers.
In case the characteristic of $K$ is zero,
let $d' = d$, $i = 1$.
For notational simplicity
we assume that $K$ contains all the ${d \over i}$th roots of unity.
Let
$s,f,s_1,f_1,c_1,\ldots,c_4,b_1,\ldots,b_4$ be indeterminates over $K$,
and $R = K[s,f,s_1,f_1,c_1,\ldots,c_4,b_1,\ldots,b_4]$.
Note that $R$ has dimension $12$.
The Mayr-Meyer ideal for $n = 1$ is the ideal in $R$
with the generators as follows:
$$
\eqalignno{
J = J(1,d) &=
(s_1-sc_1, f_1-sc_4) +
\left(c_i \left(s -fb_i^d\right)| i = 1, 2, 3, 4\right) \cr
&\hskip 2em
+ \left(fc_1 - s c_2, fc_4 - s c_3, s \left(c_3 - c_2 \right),
f \left(c_2 b_1 - c_3 b_4 \right),
f c_2 \left(b_2 - b_3 \right) \right). \cr
}
$$

\thmline{Theorem 1}
A minimal primary decomposition of $J = J(1,d)$ is as follows:
$$
\eqalignno{
J &= \left(s_1-sc_1, f_1-sc_4,c_1, c_2, c_3, c_4\right) \cr
&\hskip 1em
\bigcap_{\alpha} \left(s_1-sc_1, f_1-sc_4,c_4-c_1,c_3-c_2,c_1-c_2b_1^d,
s-fb_1^d,b_1-b_4,b_2-b_3,b_1^i-\alpha b_2^i \right) \cr
&\hskip 1em \cap
\left(s_1-sc_1, f_1-sc_4,s,f\right) \cr
&\hskip 1em \cap
\left(s_1-sc_1, f_1-sc_4,s,c_1,c_2,c_4,b_3^d,b_4\right) \cr
&\hskip 1em \cap
\left(s_1-sc_1, f_1-sc_4,s,c_1,c_4,b_3^d,b_2- b_3, c_2b_1-c_3b_4\right) \cr
&\hskip 1em \cap
\left(s_1-sc_1, f_1-sc_4,s^2,f^2,c_4(s-fb_4^d),
c_3(s-fb_3^d), sc_3-fc_4, c_3^2,c_4^2,
\right. \cr
&\hskip 4em \left.
c_1-c_4,c_2-c_3, b_2- b_3, b_1- b_4\right), \cr
}
$$
where the $\alpha$ vary over the ${d \over i}$th roots of unity in $K$.
\endb

It is easy to verify that $J = J(1,d)$ is contained in the intersection,
and that all but the last ideal
on the right-hand side of the equality are primary.
The following lemma proves that the last ideal is primary as well:

\thmline{Lemma 2}
The last ideal in the intersection in Theorem~1 is primary.
\endb

\proof
Here is a simple fact:
let $x_1, \ldots, x_n$ be variables over a ring $A$,
$S = A[x_1, \ldots, x_n]$,
and $I$ an ideal in $A$.
Then $I$ is primary (respectively, prime)
if and only if for any $f_1, \ldots, f_n \in A$,
$IS + (x_1 - f_1, \ldots, x_n - f_n)S$ is a primary (respectively, prime)
ideal in $S$.

By this fact it suffices to prove that the ideal
$$
L = \left(s^2,f^2,c_4(s-fb_4^d),
c_3(s-fb_3^d), sc_3-fc_4, c_3^2,c_4^2\right)
$$
is primary.
Note that $\sqrt L = (s,f,c_3,c_4)$ is a prime ideal.
It suffices to prove that the set of associated primes of $L$ is $\{\sqrt L\}$.
It is an easy fact
that for any $x \in R$,
$$
\Ass\left({R \over L}\right) \subseteq
\Ass\left({R \over L : x}\right) \cup
\Ass\left({R \over L+ (x)}\right).
$$
In particular,
when $x = f$,
$L + (f) = \left(s^2,f,sc_4, sc_3, c_3^2,c_4^2\right)$
is clearly primary to $\sqrt L$.
Thus it suffices to prove that $L : f$ is primary to $\sqrt L$.

We fix the monomial lexicographic ordering
$s > f > c_4 > c_3 > b_4 > b_3$.
Clearly $L : f$ contains $(s^2,f,c_4-c_3b_3^d, sc_3,c_3^2)$.
If $r \in L : f$,
then the leading term of $r$ times $f$ is contained in the ideal
of leading terms of $L$,
namely in $(s^2,f^2,sc_4, sc_3,c_3^2,c_4^2,fc_4)$,
so that
the leading term of $r$ lies in $(s^2,f,c_4, sc_3,c_3^2)$.
This proves that
$L : f = (s^2,f,c_4-c_3b_3^d, sc_3,c_3^2)$.
This ideal is clearly primary to $\sqrt L$,
which proves the lemma.
\qed

We next prove that the intersection of ideals in Theorem~1 equals $J = J(1,d)$.
Note that it suffices to prove the shortened equality:
$$
\eqalignno{
&\hskip-1em \left(c_1, c_2, c_3, c_4\right) \cr
&\hskip 1em
\bigcap_{\alpha} \left(c_4-c_1,c_3-c_2,c_1-c_2b_1^d,
s-fb_1^d,b_1-b_4,b_2-b_3,b_1^i-\alpha b_2^i \right) \cr
&\hskip 1em \cap
\left(s,f\right) \cr
&\hskip 1em \cap
\left(s,c_1,c_2,c_4,b_3^d,b_4\right) \cr
&\hskip 1em \cap
\left(s,c_1,c_4,b_3^d,b_2- b_3, c_2b_1-c_3b_4\right) \cr
&\hskip 1em \cap
\left(s^2,f^2,c_4(s-fb_4^d),
c_3(s-fb_3^d), sc_3-fc_4, c_3^2,c_4^2,
c_1-c_4,c_2-c_3, b_2- b_3, b_1- b_4\right) \cr
&\hskip-1em
=\left(c_i \left(s -fb_i^d\right),
fc_1 - s c_2, fc_4 - s c_3, s \left(c_3 - c_2 \right),
f \left(c_2 b_1 - c_3 b_4 \right),
f c_2 \left(b_2 - b_3 \right) \right). \cr
}
$$
The intersection of the first two rows equals:
$$
\eqalignno{
&\hskip-1em \left(c_1, c_2, c_3, c_4\right) 
\cap \left(c_4-c_1,c_3-c_2,c_1-c_2b_1^d,
s-fb_1^d,b_1-b_4,b_2-b_3,b_1^d-b_2^d\right) \cr
&=
\left(c_4-c_1,c_3-c_2,c_1-c_2b_1^d\right)
+
\left(c_1, c_2, c_3, c_4\right)
\cdot \left(s-fb_1^d,b_1-b_4,b_2-b_3,b_1^d-b_2^d\right) \cr
&=
J + \left(c_4-c_1,c_3-c_2,c_1-c_2b_1^d\right)
+c_2\cdot \left(b_1-b_4,b_2-b_3,b_1^d-b_2^d\right). \cr
}
$$
This intersected with the third row, namely with $(s,f)$,
equals
$$
J + (s,f)\left(c_4-c_1,c_3-c_2,c_1-c_2b_1^d\right)
+c_2(s,f)\left(b_1-b_4,b_2-b_3,b_1^d-b_2^d\right).
$$
Modulo $J$,
$$
\eqalignno{
sc_1 &\equiv fc_1b_1^d \equiv sc_2b_1^d
\equiv fc_2b_1^db_2^d
\equiv fc_2b_1^db_3^d \cr
&\equiv
fc_3b_1^{d-1}b_4b_3^d
\equiv sc_3b_1^{d-1}b_4
\equiv sc_2b_1^{d-1}b_4
\equiv fc_2b_1^{d-1}b_4b_2^d
\equiv fc_2b_1^{d-1}b_4b_3^d \cr
&\equiv fc_3b_1^{d-2}b_4^2b_3^d
\equiv \cdots \equiv fc_3b_1^0 b_4^db_3^d
\equiv sc_3b_4^d
\equiv
fc_4b_4^d
\equiv sc_4, \cr
}
$$
so that $s(c_1-c_4) \in J$.
Also it is clear that $f(c_1-c_4) , s(c_3-c_2),
s(c_1-c_2b_1^d) \in J$,
that $sc_2 \in J + (fc_2)$,
and that $fc_2(b_2-b_3) \in J$.
Thus the intersection of the ideals in the first three rows of
Theorem~1 simplifies to
$J + f\left(c_3-c_2,c_1-c_2b_1^d\right) +fc_2\left(b_1-b_4,b_1^d-b_2^d\right)$.
Furthermore,
$fc_2(b_1-b_4) \in (f(c_3-c_2)) + J$
and modulo $J$,
$f(c_1-c_2b_1^d) \equiv c_2(s-fb_1^d) \equiv
fc_2(b_2^d-b_1^d)$,
so that finally the intersection of the first three rows simplifies to
$$
J + \left(f(c_3-c_2),fc_2(b_1^d-b_2^d)\right).
$$
We intersect this with the (shortened) ideal
in the fourth row of Theorem~1,
namely with $(s,c_1,$
$c_2,c_4,b_3^d,b_4)$, to get
$$
\eqalignno{
J &+ (fc_2(b_1^d-b_2^d)) + (f(c_3-c_2)) \cap
(s,c_1,c_2,c_4,b_3^d,b_4) \cr
&=
J + (fc_2(b_1^d-b_2^d)) + f(c_3-c_2)\cdot
(s,c_1,c_2,c_4,b_3^d,b_4) \cr
&=
J + (fc_2(b_1^d-b_2^d)) + f(c_3-c_2)\cdot (c_2,b_3^d,b_4). \cr
}
$$
As modulo $J$,
$f(c_3-c_2)b_4 \equiv fc_2(b_1-b_4)$,
and $f(c_3-c_2)b_3^d \equiv fc_3b_3^d-fc_2b_2^d
\equiv sc_3-sc_2 \equiv 0$,
the intersection of the first four rows simplifies to
$$
J + fc_2\left(b_1^d-b_2^d,c_3-c_2,b_1-b_4\right).
$$
Next we intersect this with the ideal in the fifth row
(of Theorem~1)
namely with
$(s,c_1,c_4,b_3^d,b_2- b_3, c_2b_1-c_3b_4)$,
to get:
$$
\eqalignno{
J &+
(fc_2\left(b_1^d-b_2^d,c_3-c_2,b_1-b_4\right))
\cap \left(s,c_1,c_4,b_3^d,b_2- b_3,
c_2b_1-c_3b_4\right) \cr
&\hskip-1em
=J + fc_2\left(\left(b_1^d-b_2^d,c_3-c_2,b_1-b_4\right)
\cap \left(s,c_1,c_4,b_3^d,b_2- b_3,
c_2b_1-c_3b_4\right) \right)\cr
&\hskip-1em
=J + fc_2\left(\left(c_2b_1-c_3b_4\right)%
\kern-.2em+\kern-.2em\left(b_1^d-b_2^d,c_3-c_2,b_1-b_4\right)%
\kern-.2em\cap\kern-.2em\left(s,c_1,c_4,b_3^d,b_2- b_3\right)
\right)\cr
&\hskip-1em
=J + fc_2\left(\left(c_2b_1-c_3b_4\right)
+ \left(b_1^d-b_2^d,c_3-c_2,b_1-b_4\right)
\cdot \left(s,c_1,c_4,b_3^d,b_2- b_3\right)
\right)\cr
&\hskip-1em
=J + fc_2
\left(b_1^d-b_2^d,c_3-c_2,b_1-b_4\right)
\cdot \left(s,c_1,c_4,b_3^d\right). \cr
}
$$
As modulo $J$,
$$
\eqalignno{
sc_2(b_1-b_4) &\equiv fc_3b_3^d(b_1-b_4)
\equiv fb_3^d(c_3-c_2)b_1 \equiv (sc_3-fb_2^dc_2)b_1
\equiv s(c_3-c_2)b_1 \equiv 0, \cr
sfc_2(b_1^d-b_2^d) &\equiv f^2c_1b_1^d-s^2c_2 \equiv sfc_1-sfc_1 = 0, \cr
fc_1 &\equiv fc_4, \cr
}
$$
the intersection of the ideals in the first five rows simplifies to
$$
\eqalignno{
J &+ fc_2
\left(b_1^d-b_2^d,c_3-c_2,b_1-b_4\right)
\cdot \left(c_1,b_3^d\right) \cr
&=J + sc_2\left(b_1^d-b_2^d,c_3-c_2,b_1-b_4\right) \cr
&=J + sc_2\left(b_1^d-b_2^d\right). \cr
}
$$
Finally we intersect this intersection of the ideals in the first five rows
in the statement of Theorem~1 %
with the (shortened) last ideal there,
namely with
$L = (s^2,f^2,c_4(s-fb_4^d),
c_3(s-fb_3^d), sc_3-fc_4, c_3^2,c_4^2,
c_1-c_4,c_2-c_3, b_2- b_3, b_1- b_4)$, to get:
$$
J + sc_2\left(b_1^d-b_2^d\right) \cap L.
$$
It is easy to see that $L : sc_2$ contains $\sqrt{L}$.
As $sc_2$ is not in $L$,
then $L : sc_2=\sqrt{L}$,
so that the intersection of all the ideals in Theorem~1 equals
$$
\eqalignno{
J + &sc_2\left(b_1^d-b_2^d\right)\Bigl(L:sc_2\left(b_1^d-b_2^d\right)\Bigr)
=J + sc_2\left(b_1^d-b_2^d\right)\Bigl(\sqrt{L}:\left(b_1^d-b_2^d\right)\Bigr)
\cr
&=J + sc_2\left(b_1^d-b_2^d\right)\sqrt{L} \cr
&=J + sc_2\left(b_1^d-b_2^d\right)
\left(s,f,c_1,c_2,c_3,c_4,b_2-b_3,b_1-b_4\right) \cr
&=J + sc_2\left(b_1^d-b_2^d\right)
\left(f,c_2\right). \cr
}
$$
It has been proved that
$sfc_2\left(b_1^d-b_2^d\right) \in J$,
and similarly 
$sc_2^2\left(b_1^d-b_2^d\right) \in J$.
This proves that the intersection of all the listed ideals in
Theorem~1 does equal $J$.
\qed

In order to finish the proof of Theorem~1,
it remains to prove that none of the listed components is redundant.
The last component is primary to a non-minimal prime,
whereas there are no inclusion relations among the rest of the
primes.
Thus the first $d'+4$ listed components belong to minimal primes
and are not redundant.
With this it suffices to prove that $J$ has an embedded prime:

\thmline{Lemma 3}
When $n = 1$,
$c_4(s-fb_3^d)$ is in every minimal component but not in $J$.
Thus there exists an embedded component.
\endb

\proof
It has been established that
$c_4(s-fb_3^d)$ is in every minimal component.
%
Suppose that $c_4(s-fb_3^d)$ is in $J$.
Then
$$
\eqalignno{
c_4(s-fb_3^d)
&=
\sum_{i=1}^4 r_i c_i \left(s -fb_i^d\right)
+ r_5 (fc_1 - s c_2)
+ r_6 (fc_4 - s c_3)
+ r_7 s \left(c_3 - c_2 \right) \cr
&\hskip 3em
+ r_8 f \left(c_2 b_1 - c_3 b_4 \right)
+ r_9 f c_2 \left(b_2 - b_3 \right), \cr
}
$$
for some elements $r_i$ in the ring.
By the homogeneity of all elements in the two sets of variables $\{s,f\}$
and $\{c_1,c_2,c_3,c_4\}$,
without loss of generality each $r_i$ is an element of $K[b_i | i=1,2,3,4]$.
Therefore the coefficients of the $fc_i$, $sc_i$ yield the following equations:
\def\so{\hbox{\ so \ }}
$$
\eqalignno{
sc_4\cl & 1= r_4, \cr
fc_4\cl & -b_3^d = -r_4 b_4^d+ r_6, \so r_6 =b_4^d -b_3^d, \cr
fc_3\cl & 0 = r_3 b_3^d+ r_8 b_4, \so r_3 = rb_4, r_8 = - r b_3^d
\hbox{\ for some $r \in R$}, \cr
sc_3\cl & 0 = b_4^d-b_3^d -r_3 - r_7, \so r_7 = b_4^d-b_3^d- r b_4, \cr
sc_1\cl & 0 = r_1, \cr
fc_1\cl & 0 = -r_1 b_1^d + r_5, \so r_5 = 0, \cr
sc_2\cl & 0 = r_2-b_4^d+b_3^d+rb_4, \so r_2=  b_4^d-b_3^d-rb_4, \cr
fc_2\cl & 0 = -r_2b_2^d - rb_3^d b_1 + r_9(b_2 - b_3). \cr
}
$$
After expanding $r_2$ in the last equation,
$0 =-(b_4^d-b_3^d-rb_4) b_2^d- r b_3^d b_1+ r_9 \left(b_2 - b_3 \right)$,
so that $b_2^d b_3^d \in (b_1, b_4,b_2-b_3)$,
which is a contradiction.
\qed

As one embedded component has been established,
this proves the Theorem.
Thus in the case $n = 1$,
the Mayr-Meyer ideal $J(1,d)$ has $d' + 4$ minimal primes and one embedded one,
and these associated prime ideals are as follows
($\alpha$ varies over the $d'\,$th roots of unity):

\vskip 2ex
\halign{\mvrule \hskip 0.4em #\hfil \hskip 0.1em & \mvrule \relax \hfil# \hskip -0.1em \mvrule \cr
\noalign{\hrule}
associated prime ideal & \hskip0.3em height \cr
\noalign{\hrule}
\noalign{\hrule}
$\left(s_1-sc_1, f_1-sc_4,c_1, c_2, c_3, c_4\right)$ & 6 \cr
\noalign{\hrule}
$\left(s_1-sc_1, f_1-sc_4,c_4-c_1,c_3-c_2,c_1-c_2b_1^d,
s-fb_1^d,b_1-b_4,b_2-b_3,b_1-\alpha b_2\right)$ & 9 \cr
\noalign{\hrule}
$\left(s_1-sc_1, f_1-sc_4,s,f\right)$ & 4 \cr
\noalign{\hrule}
$\left(s_1-sc_1, f_1-sc_4,s,c_1,c_2,c_4,b_3,b_4\right)$ & 8 \cr
\noalign{\hrule}
$\left(s_1-sc_1, f_1-sc_4,s,c_1,c_4,b_2,b_3, c_2b_1-c_3b_4\right)$ & 8 \cr
\noalign{\hrule}
$\left(s_1-sc_1, f_1-sc_4,s,f,c_1,c_2,c_3,c_4, b_2- b_3, b_1- b_4\right)$
& 10 \cr
\noalign{\hrule}
}

\vskip 4ex

The proof of the Theorem also explicitly computes the intersection
of the first five rows of the primary decomposition,
so that:

\thmline{Proposition 4}
The intersection of all the minimal components of $J(1,d)$
equals $J + (sc_2(b_1^d - b_2^d))$.
\qed
\endb

Furthermore,
it is straightforward to compute the radical of $J(1,d)$:

\thmline{Proposition 5}
The radical of $J(1,d)$ equals
$J(1, d') + fb_3(c_3-c_2,c_2(b_1^{d'}-b_2^{d'}))$.
\endb

\proof
It is straightforward to compute the radical of each component.
Note that as in the previous proof
it suffices to compute the shortened intersection:
$$
\eqalignno{
& \left(c_1, c_2, c_3, c_4\right) \cr
&\hskip 1em
\bigcap_{\alpha} \left(c_4-c_1,c_3-c_2,c_1-c_2b_1^d,
s-fb_1^d,b_1-b_4,b_2-b_3,b_1-\alpha b_2 \right) \cr
&\hskip 1em \cap
\left(s,f\right) \cr
&\hskip 1em \cap
\left(s,c_1,c_2,c_4,b_3,b_4\right) \cr
&\hskip 1em \cap
\left(s,c_1,c_4,b_2,b_3, c_2b_1-c_3b_4\right) \cr
&=
\left(c_i \left(s -fb_i^d\right),
fc_1 - s c_2, fc_4 - s c_3, s \left(c_3 - c_2 \right),
f \left(c_2 b_1 - c_3 b_4 \right)\right) \cr
&\hskip2em + \left(f c_2 \left(b_2 - b_3 \right),
fb_3(c_3-c_2),
fb_3c_2(b_1^{d'}-b_2^{d'}) \right). \cr
}
$$
As in the proof of the Theorem,
the intersection of the first three rows equals
$J(1,d') +\left(f(c_3-c_2),fc_2(b_1^{d'}-b_2^{d'})\right)$.
Intersection with the ideal in the fourth row,
namely with $(s,c_1,c_2,c_4,b_3,b_4)$,
equals
$$
\eqalignno{
J(1,d') &+ \left(fc_2(b_1^{d'}-b_2^{d'})\right)
+ \left(f(c_3-c_2)\right) \cap (s,c_1,c_2,c_4,b_3,b_4) \cr
&=J(1,d') + \left(fc_2(b_1^{d'}-b_2^{d'})\right)
+ \left(f(c_3-c_2)\right) \cdot (s,c_1,c_2,c_4,b_3,b_4) \cr
&=J(1,d') + \left(fc_2(b_1^{d'}-b_2^{d'})\right)
+ \left(f(c_3-c_2)\right) \cdot (c_2,b_3,b_4). \cr
}
$$
When this is intersected with the ideal in the fifth row,
namely with $(s,c_1,c_4,b_2,b_3, c_2b_1-c_3b_4)$,
the resulting radical of $J(1,d)$ equals
$$
\eqalignno{
J(1&,d') + (fb_3(c_3-c_2)) \cr
&\hskip1em+\left(fc_2(b_1^{d'}-b_2^{d'}),fc_2(c_3-c_2),fb_4(c_3-c_2)\right)
\cap (s,c_1,c_4,b_2,b_3, c_2b_1-c_3b_4)\cr
&=J(1,d') + (fb_3(c_3-c_2)) \cr
&\hskip1em+fc_2\left(\left(b_1^{d'}-b_2^{d'},c_3-c_2, b_1-b_4\right)
\cap (s,c_1,c_4,b_2,b_3, c_2b_1-c_3b_4)\right)\cr
&=J(1,d') + (fb_3(c_3-c_2))
+fc_2\left(\left(b_1^{d'}-b_2^{d'},c_3-c_2, b_1-b_4\right)
\cdot (s,c_1,c_4,b_2,b_3)\right), \cr
}
$$
and by previous computations this simplifies to
$$
\eqalignno{
J(1,d') &+ (fb_3(c_3-c_2))
+fc_2b_3\left(b_1^{d'}-b_2^{d'},c_3-c_2, b_1-b_4\right) \cr
&=J(1,d') + fb_3(c_3-c_2,c_2(b_1^{d'}-b_2^{d'})). \eqed \cr
}
$$

Mayr and Meyer [MM] observed that whenever the element $s(c_4-c_1)$ of $J(1,d)$
is expressed as an $R$-linear combination of the given generators of $J(1,d)$,
at least one of the coefficients has degree at least $d$.
In fact, as the proposition below proves,
the degree of at least one of the coefficients is at least $2d-1$,
and this lower bound is achieved.
(See also the proof of Theorem showing that $sc_4 \equiv sc_1$
modulo $J(1,d)$.)
Mayr and Meyer also showed the analogues for $n \ge 1$,
with degrees of the coefficients depending on $n-1$ doubly exponentially.

Bayer, Huneke and Stillman questioned how much this doubly
exponential growth depends on the existence of embedded primes of $J(n,d)$,
or on the structure of the components.
The proposition below shows that at least for $n = 1$,
the facts that $J(1,d)$ has an embedded prime
and that the minimal components are not radical,
do not seem to be crucial for this property:

\thmline{Proposition 6}
Let $I$ be any ideal between $J(1,d)$ and its radical.
Then whenever the element $s(c_4-c_1)$
is expressed as an $R$-linear combination
of the minimal generators of $I$ which include all the given
generators of $J(1,d)$,
at least one of the coefficients has degree at least $2d-1$.
\endb

\proof
All the cases can be deduced from the case of $I$ being the
radical of $J(1,d)$.
To simplify the notation it suffices to replace $d$ by $d'$,
so that $I = J(1,d) + fb_3(c_3-c_2,c_2(b_1^d-b_2^d))$.
Write
$$
\eqalignno{
s(c_4-c_1)
&=
\sum_{i=1}^4 r_i c_i \left(s -fb_i^d\right)
+ r_5 (fc_1 - s c_2)
+ r_6 (fc_4 - s c_3)
+ r_7 s \left(c_3 - c_2 \right) \cr
&+ r_8 f \left(c_2 b_1 - c_3 b_4 \right)
+ r_9 f c_2 \left(b_2 - b_3 \right) + r_{10}fb_3(c_3-c_2)
+ r_{11}fb_3c_2(b_1^d-b_2^d) \cr
}
$$
for some elements $r_i$ in the ring.
Note that each of the explicit elements of $I$
is homogeneous in the two sets of variables $\{s,f\}$
and $\{c_1,c_2,c_3,c_4\}$.
Thus it suffices to prove that in degrees 1 in each of the two sets
of variables,
one of the coefficients has degree at least $2d-1$.
So without loss of generality each $r_i$ is an element of $K[b_i | i=1,2,3,4]$.
Therefore the coefficients of the $sc_i$, $fc_i$ yield
$$
\eqalignno{
sc_4\cl & 1 = r_4, \cr
fc_4\cl & 0 = -r_4 b_4^d + r_6, \so r_6 = b_4^d, \cr
sc_1\cl & -1 = r_1, \cr
fc_1\cl & 0 = -r_1b_1^d + r_5, \so r_5 = -b_1^d, \cr
sc_2\cl & 0 =r_2 + b_1^d- r_7, \so r_7 =r_2 + b_1^d, \cr
fc_2\cl & 0 =-r_2 b_2^d+ r_8 b_1+ r_9 \left(b_2 - b_3 \right) - r_{10}b_3
	+ r_{11}b_3(b_1^d-b_2^d), \cr
sc_3\cl & 0 =r_3 - b_4^d+ r_7, \so r_3 = b_4^d-b_1^d - r_2, \cr
fc_3\cl & 0 =-r_3 b_3^d- r_8 b_4+ r_{10}b_3. \cr
}
$$
The last equation implies that $r_8 = r'_8 b_3$,
so that the equations for the coefficients of $fc_3$ and $fc_2$
can be rewritten as
$$
\eqalignno{
r_{10} &= r_3 b_3^{d-1} + r'_8 b_4 =
(b_4^d-b_1^d) b_3^{d-1}- r_2 b_3^{d-1}+ r'_8 b_4, \cr
0 &=-r_2 b_2^d+ r'_8 b_3 b_1+ r_9 \left(b_2 - b_3 \right)
- (b_4^d-b_1^d) b_3^d+r_2 b_3^d- r'_8 b_4b_3+ r_{11}b_3(b_1^d-b_2^d), \cr
&=-r_2 (b_2^d-b_3^d)+ r'_8 b_3( b_1-b_4)
+ r_9 (b_2 - b_3)
-(b_4^d-b_1^d)b_3^d
+ r_{11}b_3(b_1^d-b_2^d). \cr
}
$$
Thus $r'_8 b_3(b_1-b_4) \in (b_2-b_3, b_3^d(b_1^d-b_4^d), b_1^d-b_2^d)$,
so that
$r'_8 \in (b_2-b_3, b_3^{d-1}{b_1^d-b_4^d \over b_1-b_4}, b_1^d-b_2^d)$.
If $r'_8 \in (b_2-b_3, b_1^d-b_2^d)$.
Then
$(b_4^d-b_1^d)b_3^d \in (b_2-b_3, b_1^d-b_2^d)$,
which is a contradiction.
Thus $r'_8$ has a multiple of $b_3^{d-1}{b_1^d-b_4^d \over b_1-b_4}$
as a summand,
so $r'_8$ has degree at least $2d-2$,
so that $r_8$ has degree at least $2d-1$.
In fact,
by setting all the free variables $r_2, r_{11}$ to zero,
the maximum degree of the coefficients $r_i$ is $2d-1$.
\qed

Note that in the proof above it is possible to have both
$r_{10} = r_{11} = 0$
and the degrees of the $r_i$ still at most $2d-1$,
with $2d-1$ attained on some $r_i$.
(Lemma 2.3 of [BS] erroneously claims that the degree of some $r_i$
is at least $2d$.)


\medskip

\vskip 4ex

\bigskip
\leftline{\bf References}
\bigskip

\bgroup
\font\eightrm=cmr8 \def\rm{\fam0\eightrm}
\font\eightit=cmti8 \def\it{\fam\itfam\eightit}
\font\eightbf=cmbx8 \def\bf{\fam\bffam\eightbf}
\font\eighttt=cmtt8 \def\tt{\fam\ttfam\eighttt}
\rm
\baselineskip=9.9pt
\parindent=3.6em

\item{[BS]}
Bayer, D.; Stillman, M.,
On the complexity of computing syzygies,
{\it J.\ Symbolic Comput.} {\bf 6}, (1988), 135-147.

\item{[GS]}
Grayson, D.; Stillman, M.,
Macaulay2. 1996.
A system for computation in algebraic geometry and commutative algebra,
available via anonymous {\tt ftp} from {\tt math.uiuc.edu}.

\item{[H]}
Herrmann, G.,
Die Frage der endlich vielen Schritte in der Theorie der Polynomideale,
{\it Math.\ Ann.} {\bf 95}, (1926), 736-788.

\item{[K]}
Koh, J.,
Ideals generated by quadrics exhibiting double exponential degrees,
{\it J.\ Algebra} {\bf 200}, (1998), 225-245.

\item{[MM]}
Mayr, E.; Meyer, A.,
The complexity of the word problems for commutative semigroups
and polynomial ideals,
{\it Adv.\ Math.} {\bf 46}, (1982), 305-329.

\item{[S]}
Swanson, I.,
Primary decomposition of the Mayr-Meyer ideals,
in preparation, 2001.

\egroup

\end

First assume that $I = J(1,d) +(sc_2(b_1^d-b_2^d))$,
i.e., the intersection of the minimal components of $J(1,d)$. 
Write
$$
\eqalignno{
s(c_4-c_1)
&=
\sum_{i=1}^4 r_i c_i \left(s -fb_i^d\right)
+ r_5 (fc_1 - s c_2)
+ r_6 (fc_4 - s c_3)
+ r_7 s \left(c_3 - c_2 \right) \cr
&+ r_8 f \left(c_2 b_1 - c_3 b_4 \right)
+ r_9 f c_2 \left(b_2 - b_3 \right) + r_{10}sc_2(b_1^d-b_2^d), \cr
}
$$
The coefficients of $c_2$ and $c_3$ give
$$
\eqalignno{
0 &=
r_2 \left(s -fb_2^d\right)+s b_1^d- r_7 s + r_8 f b_1
+ r_9 f \left(b_2 - b_3 \right) + r_{10}s(b_1^d-b_2^d), \cr
0 &=
r_3 \left(s -fb_3^d\right)
- s b_4^d+ r_7 s - r_8 f  b_4, \cr
}
$$
from which the coefficients of $s$ and $f$ yield
$$
\eqalignno{
0 &=r_2 + b_1^d- r_7+r_{10}(b_1^d-b_2^d), \cr
0 &=-r_2 b_2^d+ r_8 b_1+ r_9 \left(b_2 - b_3 \right), \cr
0 &=r_3 - b_4^d+ r_7, \cr
0 &=-r_3 b_3^d- r_8 b_4. \cr
}
$$
The last equation implies that $r_8 = r b_3^d$, $r_3 = - r b_4$
for some $r \in R$.
Thus the equations above imply
$$
\eqalignno{
r_7 &=b_4^d + rb_4, \cr
r_2 &= r b_4 + b_4^d-b_1^d-r_{10}(b_1^d-b_2^d), \cr
\left(r b_4 + b_4^d-b_1^d-r_{10}(b_1^d-b_2^d)\right) b_2^d
&=r b_3^d b_1+ r_9 \left(b_2 - b_3 \right) \cr
&= r b_2^d b_1+ (r_9  + rb_1 {b_3^d - b_2^d \over b_2 - b_3})
(b_2 - b_3). \cr
}
$$
Thus there exists $p \in R$ such that
$$
\eqalignno{
p (b_2 - b_3) &= r (b_4 - b_1)+ b_4^d-b_1^d-r_{10}(b_1^d-b_2^d), \cr
p b_2^d 
&= r_9  + rb_1 {b_3^d - b_2^d \over b_2 - b_3}. \cr
}
$$
If $r_9 = 0$,
then $r \in (b_2^d)$,
hence $b_4^d-b_1^d \in (b_2 - b_3, b_2^d, b_1^d-b_2^d)$,
which is a contradiction.
So necessarily $r_9 \not = 0$.
The first equation above implies that $r$ is non-zero
and of degree at least $d-1$,
and then the second equation implies that the degree of $r_9$ is at 
least $2d-1$.

In fact,
setting $r$ to be $-{b_4^d-b_1^d \over b_4 - b_1}$
and all other free variables to $0$,
the maximum degree of $2d-1$ for $r_i$ is achieved.
As $r_{10} = 0$,
this also proves the case when $I = J(1,d)$.
\qed

\prop
Let $I$ be the radical of $J(1,d)$.
Then whenever the element $s(c_4-c_1)$
is expressed as an $R$-linear combination of the given generators of $I$,
at least one of the coefficients has degree at least $2d-1$.
\endb

\proof
To simplify the notation it suffices to replace $d$ by $d/i$,
so that the radical equals $J(1,d) + fb_3(c_3-c_2,c_2(b_1^d-b_2^d))$.
Write
$$
\eqalignno{
s(c_4-c_1)
&=
\sum_{i=1}^4 r_i c_i \left(s -fb_i^d\right)
+ r_5 (fc_1 - s c_2)
+ r_6 (fc_4 - s c_3)
+ r_7 s \left(c_3 - c_2 \right) \cr
&+ r_8 f \left(c_2 b_1 - c_3 b_4 \right)
+ r_9 f c_2 \left(b_2 - b_3 \right) + r_{10}fb_3(c_3-c_2)
+ r_{11}fb_3c_2(b_1^d-b_2^d), \cr
}
$$
for some elements $r_i$ in the ring.
As in the previous proof
it suffices to consider the case when
each $r_i$ is an element of $K[b_i | i=1,2,3,4]$.
Isolating the coefficients of $c_4$ and $c_1$
as in the previous proof gives
$r_4 = 1$ and $r_6 =b_4^d$, $r_1 = -1$, $r_5 = -b_1^d$, and
$$
\eqalignno{
0 &=
\sum_{i=2}^3 r_i c_i \left(s -fb_i^d\right)
+s c_2 b_1^d- s c_3b_4^d
+ r_7 s \left(c_3 - c_2 \right) \cr
&\hskip 3em
+ r_8 f \left(c_2 b_1 - c_3 b_4 \right)
+ r_9 f c_2 \left(b_2 - b_3 \right) + r_{10}fb_3(c_3-c_2)
+ r_{11}fb_3c_2(b_1^d-b_2^d). \cr
}